%% file: main.tex
\documentclass{ifacconf}

\usepackage{graphicx}      
\graphicspath{{figurescodedata/}}
\usepackage{natbib}        
\usepackage{amsfonts,mathtools}       %
\usepackage{nicefrac} 

\usepackage{microtype}      
\usepackage{xcolor,color}         
\usepackage{amsmath,bm}
\usepackage{amssymb}

\allowdisplaybreaks

\usepackage{subcaption}
\usepackage{booktabs} 
\usepackage{enumerate}
\usepackage{enumitem}
\usepackage{setspace}

\usepackage{algorithm}
\usepackage{algorithmic}
\usepackage[algo2e,ruled]{algorithm2e}

\DeclareMathOperator*{\argmin}{arg\,min}

\DeclareMathOperator*{\E}{\mathbb{E}}

\theoremstyle{plain}

\newtheorem{theorem}{Theorem}[section]

\newtheorem{lemma}[theorem]{Lemma}
\newtheorem{corollary}[theorem]{Corollary}
\newtheorem{example}{Example}
\newtheorem{assumption}[theorem]{Assumption}

\theoremstyle{definition}
\newtheorem{definition}[theorem]{Definition}

\theoremstyle{plain}

\definecolor{dark}{rgb}{0.35, 0.15, 0.13}

\newcommand{\R}{\mathbb R}

\newcommand{\U}{\mathbb U}
\newcommand{\W}{\mathbb W}

\newcommand{\F}{\mathcal F}

\newcommand{\X}{\mathbb X}
\newcommand{\Pb}{\mathbb P}

\newcommand{\Sb}{\mathbb S}

\usepackage{dsfont}
\newcommand{\one}{\mathds{1}}

\newcommand{\nbf}{\noindent\textbf}

\begin{document}
\begin{frontmatter}

\title{Non-asymptotic System Identification  for   Linear Systems with Nonlinear  Policies} 

\thanks[footnoteinfo]{This work is supported by NSF AI institute: 2112085 and ONR YIP: N00014-19-1-2217.}

\author[1,2]{Yingying Li} 
\author[1]{Tianpeng Zhang} 
\author[3]{Subhro Das}
\author[2]{Jeff Shamma}
\author[1]{Na Li}

\address[1]{Harvard University}
\address[2]{University of Illinois at Urbana-Champaign}
\address[3]{MIT-IBM Watson AI Lab, IBM Research}

\begin{abstract}                
This paper considers a single-trajectory system identification problem for linear systems under general nonlinear and/or time-varying policies with i.i.d. random excitation noises. The problem is motivated by safe learning-based control for constrained linear systems, where the safe policies during the learning process are usually nonlinear and time-varying for satisfying the state and input constraints. In this paper, we provide a non-asymptotic error bound for least square estimation when the data trajectory is generated by any nonlinear and/or time-varying  policies as long as the generated state and action trajectories are bounded. This significantly generalizes the existing non-asymptotic guarantees for linear system identification, which usually consider i.i.d. random inputs or linear policies. Interestingly, our error bound is consistent with that for linear policies with respect to the dependence on the trajectory length, system dimensions, and excitation levels. Lastly, we demonstrate the applications of our results by safe learning with robust model predictive control and provide numerical analysis.
\end{abstract}


\end{frontmatter}
\input{sec1_introduction}
\input{sec2_problem}

\input{sec3_algorithm}

\input{sec4_theory_proof}
\input{sec5_numerical}

\input{sec6_conclusion}

\bibliography{citation4safeID}             

\end{document}

%% file: sec1_introduction.tex
\section{Introduction}


This paper considers a system identification problem of a linear system $x_{t+1}=A_* x_t+B_* u_t+w_t$ under a single trajectory of data $\{ x_t, u_t\}_{t=0}^T$ generated by potentially nonlinear, time-varying, and/or history-dependent policies:
\begin{align}\label{equ: sec1 ut phit def}
	 u_t =\pi_t(x_t, \{x_s, w_s, \eta_s\}_{s=0}^{t-1})+\eta_t,
\end{align}
where  $\eta_t$ is included to provide excitation for the system exploration. The disturbances $w_t$ and $\eta_t$ are i.i.d. and bounded. We do not impose any structural assumptions on the policies $\pi_t$  except that the policies generate bounded state and action trajectories. 
We adopt the least square estimation. Our goal is to provide  non-asymptotic bounds for the estimation errors. 

Though the single-trajectory identification of linear systems  has been well-studied when the data are generated from, e.g., i.i.d. random control inputs \citep{simchowitz2018learning,dean2019sample}, and  linear policies \citep{dean2018regret,dean2019safely}, the identification with data generated from \textit{nonlinear} and \textit{time-varying} policies is relatively less explored. 

The motivation for considering nonlinear and time-varying policies for linear systems with bounded states and actions is from  safe learning of linear systems with state and action constraints \citep{lorenzen2019robust,kohler2019linear,rawlings2009model}, which enjoys wide applications in safety-critical applications \citep{fisac2018general}. Many control designs with robust constraint satisfaction   despite model uncertainties will generate nonlinear policies even for linear systems, e.g., robust model predictive control (RMPC) \citep{rawlings2009model}, the controllers based on control-barrier functions (CBF) \citep{xu2018constrained}, etc. Further, the policies can even be time-varying when the objective is time-varying, e.g., in the tracking problems \citep{limon2010robust}, and/or when the model uncertainties are adaptively updated, e.g., robust adaptive MPC \citep{lorenzen2019robust,kohler2019linear}. Therefore, to better understand the learning performance of safe controllers in constrained linear systems, it is crucial to study more general policy classes such as \eqref{equ: sec1 ut phit def} and analyze  the corresponding non-asymptotic estimation performance.

Since our closed-loop system is nonlinear under the policy \eqref{equ: sec1 ut phit def}, our problem is also related to the non-asymptotic  identification  of  nonlinear systems (see, e.g., \citet{ziemann2022learning,foster2020learning}), especially  generalized linear systems, which has received a lot of attention recently due to its connection with neural networks (see, e.g., \citet{oymak19stochastic,mania2022active,sattar2022non}). 
A majority of papers in this area focus on a different setting from ours, i.e., $x_{t+1}=\sigma(A_* x_t)+w_t$, where the nonlinear component $\sigma(\cdot)$ is outside the unknown linear component \citep{sattar2022non,foster2020learning}. Further, these papers usually require the closed-loop system to be exponentially stable, e.g., \citep{sattar2022non,foster2020learning}, which may not be satisfied by the safe control with model uncertainties.\footnote{For example, though tube-based robust MPC enjoys exponential stability when  $A_*, B_*$ are known (Prop. 3.15 in \citet{rawlings2009model}), it only has an asymptotic stability guarantee when $A_*, B_*$ are unknown (Prop. 3.20 in \citet{rawlings2009model}).} In contrast, \citet{mania2022active} focus on a system  $x_{t+1}=A_*\phi(x_t, \eta_t)+w_t$ with unknown $A_*$ and known nonlinear features $\phi(x_t, \eta_t)$, which is more general than our closed-loop system with a time-invariant and memoryless policy, i.e., $x_{t+1}=A_* x_t + B_*( \pi(x_t)+\eta_t)+w_t$. Interestingly, \citet{mania2022active} argue that, after involving a known nonlinear component $\phi$,   i.i.d. random excitation $\eta_t$ is \textit{not} enough to learn the unknown parameters $A_*$ efficiently, which is in stark contrast to closed-loop linear systems.    Thus,  \citet{mania2022active}  propose an excitation generation algorithm to obtain non-asymptotic estimation guarantees. This gives rise to some interesting questions below.

\textit{Questions: is i.i.d. random excitation enough  for linear system identification with general nonlinear and/or time-varying control policies? How much difference will the theoretical guarantee be from that of linear policies?}

\nbf{Contributions. } Our major contribution is  showing that i.i.d. excitations and the bounded state-action assumption are enough for  our system identification problem by providing a non-asymptotic estimation error bound for least square estimation of linear systems under general nonlinear and/or time-varying policies \eqref{equ: sec1 ut phit def}. 
Further, our estimation error bound scales as $\tilde O\left(\frac{\sqrt{m+n}}{\sigma_{\eta}\sqrt T}\right)$ under proper conditions, where is the same order as that of the linear policies in \citet{dean2019safely} with respect to the state dimension $n$, control input dimension $m$, the minimal eigenvalue of the covariance matrix of excitation $\eta_t$, and the length of the trajectory $T$. This indicates that, for linear systems, allowing more  general data-generation policies will not degrade the learning performance  compared with the linear policies with respect to the trajectory length, system dimensions, and excitation levels,  as long as the states and actions are bounded.
Lastly, to demonstrate the applications of this result, we consider RMPC as an example and provide its estimation error bound. We also conduct numerical experiments for safe learning with RMPC  to complement our theoretical analysis. 









\nbf{Related works. } We will review the related literature on system identification and  constrained control below.

\textit{System identification. } 
System identification  enjoys a long history of research (see e.g., \citet{boyd1986necessary,fogel1982value}). This work is mostly related to the non-asymptotic analysis of least-square estimation for linear dynamical systems, including linear system identification with linear policies and unbounded disturbances \citep{simchowitz2018learning,dean2018regret,dean2019sample}, linear system identification with bounded disturbances and robust constraint satisfaction by linear policies \citep{dean2019safely}, etc. Here, we also consider  bounded disturbances and robust constraint satisfaction but extends the results to general nonlinear, time-varying, and/or history-dependent policies. 


There is also a growing interest in the identification of nonlinear systems, e.g., bilinear systems \citep{sattar2022finite}, generalized linear systems \citep{mania2022active,foster2020learning,oymak19stochastic,sattar2022non}, or general nonlinear systems \citep{ziemann2022learning,foster2020learning}. In this paper, we consider a special closed-loop nonlinear system motivated by safe learning of constrained linear systems and leverage the structures of our problem to relax the assumptions such as exponential stability in \citep{foster2020learning,sattar2022non} and show that i.i.d. random excitation is enough for estimation in our case, though it is not enough in the general case in \citep{mania2022active}. 


It is worth mentioning a line of work that aims to design efficient data generation policies to improve the non-asymptotic estimation guarantee, e.g., \citep{zhao2022adaptive,mania2022active}, which can be viewed as an orthogonal direction of this work since this paper aims to provide estimation guarantee for as general policies as possible.

Another popular  identification method is set membership \citep{bai1998convergence,fogel1982value}. In the literature on safe learning for constrained linear systems, both set membership and least square estimation have been adopted  \citep{lorenzen2019robust}.

Lastly, there are    non-asymptotic analysis of linear system identification  with output feedback, e.g., \citep{mhammedi2020learning, sarkar2019near,oymak2019non}.




\textit{Robust  control with constraints. }  Popular methods for  robust control with constraint satisfaction include, e.g., RMPC \citep{lorenzen2019robust,kohler2019linear,rawlings2009model}, control barrier functions (CBF) \citep{salehi2022learning,xu2018constrained,taylor2020control,lopez2020robust}, safety certification \citep{wabersich2018linear,fisac2018general}, system level synthesis \citep{dean2019safely}, disturbance-action policies \citep{li2021online,li2021safe} etc. Among them, RMPC, CBF-based methods, and control with safety certification all generate potentially nonlinear policies, and system-level synthesis will generate linear policies depending on history. Notice that our policy form \eqref{equ: sec1 ut phit def} includes all these policies. 

Recent years have witnessed great interest in   safe adaptive learning for robust control with constraints, e.g., \citep{lorenzen2019robust,kohler2019linear,dean2019safely,fisac2018general}. The non-asymptotic regret analysis for this problem has also attracted  growing attention.  For example, \citet{wabersich2020performance} adopts a Thompson sampling approach, but the computation of posterior distributions can be demanding.
To ease the computation issue, $\epsilon$-greedy or certainty-equivalence type of approaches are commonly adopted for learning-based control without constraints \citep{dean2018regret}. Recently, \citet{dogan2021regret}  explored this direction and provided a regret analysis for RMPC. However, they utilize i.i.d. random control inputs for exploration and switch to nonlinear RMPC policies with some  probability for exploitation. In this paper, instead of abruptly switching policies, we provide a general estimation error bound generated by the summation of the nonlinear  RMPC policies and the i.i.d. excitation noises, which lay a foundation for future regret analysis of this control design.






\nbf{Notations.}
For a matrix $\Sigma\in \R^{n\times n}$, let $\sigma_{\min}(\Sigma)$ and $\sigma_{\max}(\Sigma)$  denote the minimal and the maximum singular value, respectively. Let $I_n$ denote the identity matrix in $\R^{n\times n}$. For a random vector $x\in\R^n$, let $\text{cov}(x)$ denote its covariance matrix. Let $\one_A$ denote an indicator function on set $A$. For $a, b\in \R$, we write $a\lesssim b$ if $a\leq cb$ for some absolute constant $c>0$ and $A\preceq B$ if matrix $B-A$ is positive semidefinite. Define $\mathbb X \oplus \mathbb Y\!=\!\{x\!+\!y: x\in \mathbb X, y\in \mathbb \! Y\! \}$ and the same for $\ominus$.

%% file: sec2_problem.tex
\vspace{-2pt}
\section{Problem formulation}\label{sec: problem}

In this paper, we consider a linear dynamical system with unknown system parameters $(A_*, B_*)$ as described below.
\begin{align}\label{equ: LTI system}
	x_{t+1}=A_* x_t + B_* u_t +w_t,
\end{align}
where $x_t \in \R^n$ and $u_t\in \R^m$.
For notational simplicity, we let $\theta_*=(A_*, B_*)$ and denote $z_t=(x_t^\top, u_t^\top)^\top$, then the system \eqref{equ: LTI system} can be written as $x_{t+1}=\theta_* z_t +w_t$.

We adopt the least square estimator (LSE) as defined below to estimate the unknown system parameters.
\begin{align}\label{equ: LSE}
	(\hat A,  \hat B)= \argmin_{A, B}\sum_{s=1}^T\|x_{s}-Ax_{s-1}-Bu_{s-1}\|_2^2.
\end{align}
For notational simplicity, we denote $\hat \theta=(\hat A, \hat B)$ and $\theta=(A, B)$. 

Our goal is to provide a non-asymptotic analysis for the  errors of LSE given a finite trajectory of states and actions generated by general policy forms described below.
\begin{align}\label{equ: ut pit}
	u_t = \tilde u_t +\eta_t, \quad \tilde u_t = \pi_t(x_t, \{x_s, u_s, \eta_s\}_{s=0}^{t-1}),
\end{align}
where $\eta_t$ is a random disturbance to provide  excitation to the control inputs, and the nominal control input $\tilde u_t$ can be generated by general policies $\pi_t(x_t, \{x_s, u_s, \eta_s\}_{s=0}^{t-1})$, which can be nonlinear, time-varying, and/or depend on the history. 
The major contribution of this paper is to provide an estimation error bound for such general policies. The only requirement on the policies is that the policy sequence will induce bounded states and control inputs.


\begin{assumption}\label{ass: bounded state action}
	The states and actions trajectories generated by the closed-loop system induced by \eqref{equ: LTI system} and \eqref{equ: ut pit} are bounded almost surely (a.s.), i.e., there exists $b_z\geq 0$ such that $\max_{t\geq 0}\|z_t\|_2=\max_{t\geq 0}\sqrt{\|x_t\|_2^2+ \|u_t\|_2^2} \leq b_z$ a.s..
\end{assumption}

This problem is motivated by  safe learning for constrained linear systems. In particular, consider state and input constraints:
\begin{align}\label{equ: constraints xt ut}
	x_t \in \X, \quad u_t \in \U,
\end{align}
where $\X, \U$ are bounded. A common question in safe adaptive learning for control is to learn $(A_*, B_*)$ without violating the constraints. A lot of control policies have been proposed with constraint satisfaction guarantees despite uncertainties in the system, thus satisfying our Assumption \ref{ass: bounded state action}. We list some safe control designs  below as examples.

\begin{example}[RMPC]
	 RMPC is commonly used to satisfy the state and action constraints, e.g., \eqref{equ: constraints xt ut}, in the presence of uncertainties in the system, e.g., disturbances $w_t$, excitation noises $\eta_t$,  uncertainties in $A_*, B_*$, etc \citep{lorenzen2019robust,kohler2019linear,rawlings2009model}. Hence, the RMPC controller with random excitation, denoted by $u_t=\pi_{\text{RMPC}}(x_t) +\eta_t$,  can satisfy Assumption \ref{ass: bounded state action} under proper conditions. Note that the RMPC controller $\pi_{\text{RMPC}}(x)$ can be nonlinear even for linear systems. In particular, $\pi_{\text{RMPC}}(x)$ is shown to be piecewise affine in $x$ for linear systems if the constraints \eqref{equ: constraints xt ut} are polytopes.
	 
	 RMPC controller can also be time-varying, e.g., when tracking a time-varying target \citep{limon2010robust}, and/or when adaptively updating the policy with improved model estimations  \citep{lorenzen2019robust,kohler2019linear}.
	 

	
\end{example}

\begin{example}[Control barrier function (CBF)] CBF is also a popular method to satisfy state and action constraints despite uncertainties and excitation noises in the system  \citep{taylor2020control,lopez2020robust}. Similar to RMPC, CBF controllers can  also be nonlinear even for linear systems and are piecewise affine for linear systems with polytopic constraints.  Hence, CBF controllers can also satisfy our Assumption \ref{ass: bounded state action}. 
	%
	%
	%
	%
	%
	%
\end{example}

\begin{example}[System-level-synthesis (SLS)]
	SLS has also been adopted to ensure constraint satisfaction under model uncertainties in linear systems \citep{dean2019safely}. Notice that the SLS controllers depend on the history states even for state feedback, which motivates us to allow policies with memory in Assumption \ref{ass: bounded state action}. In \citep{dean2019safely}, a non-asymptotic system identification error bound has been proposed for a time-invariant SLS policy. This work can complement the result in \citep{dean2019safely} by allowing time-varying SLS policies. 
	
\end{example}

%

\begin{example}[Safety certification] Safety certification has also been adopted in safe learning-based control in combination with other approaches without safety guarantees \citep{fisac2018general,wabersich2018linear}. Such algorithm design adopts classical learning approaches in the interior of the safe region and switches to a safe policy under certain criteria, e.g., on/near the boundary of the safe region \citep{fisac2018general}. Such a switching-based algorithm design naturally generates time-varying and possibly nonlinear policies, which are included by \eqref{equ: ut pit} and satisfy our Assumption \ref{ass: bounded state action}.
	

\end{example}

%

In addition, we introduce some assumptions on $w_t, \eta_t$ below.

%
%
%
\begin{assumption}[Properties of $w_t$]\label{ass: w Sigma}
	The process noise $w_t$ is i.i.d., zero mean, and bounded by $w_t \in \W =\{w: \|w\|_2\leq w_{\max}\}$. Further, the minimum eigenvalue of $\text{cov}(w_t)$ is lower bounded by $\sigma_w^2>0$. 
	

\end{assumption}




\begin{assumption}[Requirements on $\eta_t$]\label{ass: eta}
	The excitation disturbance $\eta_t$ is i.i.d., zero mean, and bounded by $\|\eta_t\|_2\leq \eta_{\max}\}$. Further, the minimum eigenvalue of $\text{cov}(\eta_t)$ is lower bounded by $\sigma_{\eta}^2>0$.
\end{assumption}
Distributions that satisfy Assumptions \ref{ass: w Sigma} and \ref{ass: eta} include, e.g., truncated Gaussian, uniform distributions on $l_2$ sphere or $l_2$ ball, etc.  The assumptions of i.i.d., zero mean and positive definite covariance matrices are commonly imposed in the linear system identification literature for non-asymptotic analysis \citep{dean2019sample,dean2019safely,simchowitz2018learning}. As for the bounded disturbances and noises, they are necessary for robust   satisfaction of bounded constraints and are thus commonly assumed in the literature of robust control with constraints \citep{rawlings2009model,lorenzen2019robust,kohler2019linear}. It may be possible to relax the boundedness assumption to Gaussian noises by considering chance constraints as in \cite{oldewurtel2008tractable}, but the verification of this relaxation is left for future work.


%
%

%% file: sec3_algorithm.tex
\section{System identification error bound}
\vspace{-3pt}

%
%
%
%
%
%

In this section, we discuss the non-asymptotic estimation error bound for least-square estimation for system \eqref{equ: LTI system} and generic nonlinear and time-varying policies \eqref{equ: ut pit}.

Before the main theorem, we introduce the notion of regularized disturbance $\eta_t/\sigma_{\eta}$ whose covariance satisfies $\text{cov}(\eta_t/\sigma_{\eta})\succeq I_m$ and whose norm satisfies $\|\eta_t/\sigma_{\eta}\|_2\leq \eta_{\max}/\sigma_{\eta}\eqqcolon \bar \eta$. We argue that the parameter pair $(\sigma_{\eta}, \bar \eta)$ is more suitable to  describe the distribution of $\eta$ than  $(\sigma_{\eta},  \eta_{\max})$ because, after rescaling the excitation disturbances to, e.g., $2\eta_t$, both $\sigma_{\eta}$ and  $\eta_{\max}$ will  change, but $\bar \eta$ will remain the same. 
Similarly, we define $\bar w\coloneqq w_{\max}/\sigma_w$. In the following, we  adopt $\sigma_{\eta}, \bar \eta, \sigma_{w}, \bar w$ for theoretical analysis and discussions. This does not cause any loss of generality. 




\begin{theorem}[Estimation error bound]\label{thm: general estimation error bdd}
	For any $0<\delta<1/3$, when the trajectory length satisfies 
	$$T\gtrsim (m\!+\!n)\max(\bar w^4\!, \bar \eta^4\!)\log(\frac{b_z}{\delta}\text{poly}_1(\bar w, \!\bar \eta, \sigma_w^{-1}, \sigma_{\eta}^{-1})) ),$$
	with probability at least $1-3\delta$, we have
	\begin{align*}
		\| \hat \theta -\theta_*\|_2  &\lesssim \frac{b_z\sqrt{m+n}}{\sqrt T \sigma_{\eta}} \text{poly}_2(\bar w, \bar \eta, \sigma_w, \sigma_{{\eta}})\\
	& \quad \cdot\sqrt{\log(\frac{b_z}{\delta})\!+\!\log(\text{poly}_1(\bar w, \!\bar \eta, \sigma_w^{-1}, \sigma_{\eta}^{-1}) )}
	\end{align*}
	where $\text{poly}_1(\bar w, \!\bar \eta, \sigma_w^{-1}, \sigma_{\eta}^{-1}) \!=\!\max(\frac{\bar w}{\sigma_w}, \frac{\bar \eta}{\sigma_{\eta}},\frac{\bar w\bar \eta}{\sigma_w\sigma_{\eta}} )\max(\bar w, \bar \eta)$ and $\text{poly}_2(\bar w, \bar \eta,\! \sigma_w, \!\sigma_{{\eta}})\!=\bar w \max(\bar w^2,\! \bar \eta^2) \max(\bar w\sigma_{\eta}, \bar \eta \sigma_w, \bar w \bar \eta)$.
\end{theorem}

Firstly, the dependence of the estimation error bound above with respect to the trajectory length and the dimensions of the system is $O\left(\sqrt{m+n}/{\sqrt T}\right)$, which is consistent with linear system identification error bound under linear policies in \citep{dean2019safely}. Further, for small enough $\sigma_{\eta}$,\footnote{\citet{dean2019safely} also assumes $\sigma_{\eta}<\sigma_w$ when analyzing their estimation error bound.} our bound depends linearly on $\tilde O(1/\sigma_{\eta})$, which is also consistent with the  bound in \citep{dean2019safely}. Further, as the process noise level $\sigma_w$ goes to infinity, the estimation error bound increases, which is also the case in the study of linear policies \citep{dean2019safely}. In summary, though we allow general nonlinear and time-varying policies to generate the data, our estimation error bounds are similar to the bound generated by linear policies. This suggests that  general data-generation policies will not significantly  degrade the estimation quality in our problem with respect to $n, m, T, \sigma_{\eta}, \sigma_w$ in comparison with the linear data-generation policies.

Besides, our estimation error bound increases with the state and action bound $b_z$. One intuitive explanation is that since our bound holds for all policies satisfying the bound $b_z$,  a larger $b_z$ includes more admissible policies, thus potentially including some  policies that generate worse estimation.

Lastly, our estimation error  increases with $\bar \eta$ and $\bar w$. This can be intuitively explained by the following:  larger $\bar \eta$ and $\bar w$ suggests more concentrated distributions in certain sense,\footnote{For example, consider the following regularized distribution: $\Pb(X=-\bar \eta)= \Pb(X=\bar \eta)=p$ and $\Pb(X=0)=1-2p$. The variance is $p\bar \eta^2= 1$, so $p=1/\bar \eta^2$. Hence, a larger $\bar  \eta$ leads to a smaller anti-concentration probability $\Pb(|X|\geq \epsilon)$ for $0<\epsilon < \bar \eta$.}  but active exploration of the unknown system calls for less concentrated disturbances, so larger $\bar \eta$ and $\bar w$ tend to provide worse estimation quality. 

\subsection{Proof of Theorem \ref{thm: general estimation error bdd}.}
The proof relies on the  block martingale small-ball (BMSB) condition introduced in \cite{simchowitz2018learning}, which is stated below for completeness.
\begin{definition}[BMSB in \citep{simchowitz2018learning} ]\label{def: bmsb}
	Let $\{\F_t\}_{t\geq 1}$ denote a filtration and 
	let  $\{Z_t\}_{t\geq 1}$ be an $\{\F_t\}_{t\geq 1}$-adapted random process taking values in $\R^d$.
	We say that  $\{Z_t\}_{t\geq 1}$ satisfies the $(k, \Gamma_{sb},p)$-block martingale small-ball (BMSB) condition for a positive integer $k$, a positive definite matrix $\Gamma_{sb}\succ 0$, and $0\leq p \leq 1$,  if the following condition holds: for  any fixed $\lambda \in\R^d$ such that $\|\lambda\|_2=1$, the process $\{Z_t\}_{t\geq 1}$  satisfies $\frac{1}{k}\sum_{i=1}^k \Pb(|\lambda^\top Z_{t+i}|\geq \sqrt{\lambda^\top\Gamma_{sb} \lambda}\mid \F_t)\geq p$ almost surely for any $t\geq 1$.
	
	%
\end{definition}

The major component of the proof is to show that the trajectory $\{z_t\}_{t\geq 0}$ satisfies the BMSB condition for general nonlinear and time-varying policies \eqref{equ: ut pit} as long as the trajectory $\{z_t\}_{t\geq 0}$ is bounded (Assumption \ref{ass: bounded state action}). By leveraging the boundedness assumption, this  result significantly relaxes the assumptions/conditions on the control policies in the literature \citep{dean2019safely,dean2018regret}.

\begin{lemma}[Verification of BMSB condition]\label{lem: verify bmsb strictly safe perturbed policy}
	Define filtration $\F_t=\{w_0, \dots, w_{t-1}, \eta_0, \dots,\eta_t\}$. Under the conditions in Theorem \ref{thm: general estimation error bdd}, 	\begin{equation*}
		\{z_t\}_{t\geq 0} \ \text{ satisfies the $(1, s_z^2 I_{n+m}, p_z)$-BMSB condition,}
	\end{equation*}
	where $p_z=\min(p_w, p_{\eta})$, $s_z=\min(s_w/4, \frac{\sqrt 3}{2}s_\eta, \frac{s_ws_\eta}{4b_z} )$, $s_w= \frac{\sigma_w}{4\bar w}$, $p_w=\frac{1}{4\bar w^2}$, $s_\eta= \frac{\sigma_\eta}{4\bar \eta}$, $p_\eta=\frac{1}{4\bar \eta^2}$.
	%
	%
	%
	%
	%
\end{lemma}
The proof is deferred to Section \ref{subsec: verify bmsb}. 


With Lemma \ref{lem: verify bmsb strictly safe perturbed policy}, we are ready for the proof of Theorem \ref{thm: general estimation error bdd}, which leverages a general least square estimation error bound in \cite{simchowitz2018learning} for time series, which is included below for completeness.
\begin{theorem}[\cite{simchowitz2018learning}]\label{thm: theorem 2.4 general estimation error}
	Fix $\epsilon\in (0,1)$, $\delta \in (0, 1/3)$, $T\geq 1$, and $0\prec\Gamma_{sb}\preceq\bar \Gamma$. Consider a random process $\{Z_t, Y_t\}_{t\geq 1}\in (\R^d\times \R^n)^T$ and a filtration $\{\F_t\}_{t\geq 1}$. 
	Suppose the following conditions hold,
	\begin{enumerate}
		\item $Y_t=\theta_* Z_t +\beta_t$, where $\beta_t\mid \F_t$ is $\sigma_{sub}^2$-sub-Gaussian and has zero mean,
		\item $\{Z_t\}_{t\geq 1}$ is an $\{\F_t\}_{t\geq 1}$-adapted random process satisfying the $(k, \Gamma_{sb},p)$-BMSB condition,
		\item $\Pb(\sum_{t=1}^T Z_t Z_t^\top \not \preceq T \bar \Gamma)\leq \delta$.
	\end{enumerate}
	If the trajectory length satisfies
	$$
	T\geq T_0=\frac{10k}{p^2}\left(\log(\frac{1}{\delta})+2d\log(10/p)+\log \det(\bar \Gamma \Gamma_{sb}^{-1})\right),
	$$
	then 	the estimation error  of the least square estimator, defined by $\hat \theta\in \argmin_{\theta}\sum_{t=1}^T \|Y_t-\theta Z_t\|^2$, satisfies
	\begin{align*}
		&	\|\hat \theta-\theta_*\|_2 \\
		& \leq \frac{90\sigma_{sub}}{p}\sqrt{\frac{n+d\log(10/p)+\log\det(\bar \Gamma \Gamma_{sb}^{-1})+\log(1/\delta)}{T \sigma_{\min}(\Gamma_{sb})}}
	\end{align*}
	with probability at least $1-3\delta$.
\end{theorem}

Now, we  will prove Theorem \ref{thm: general estimation error bdd} by verifying the conditions in Theorem \ref{thm: theorem 2.4 general estimation error}.
Condition 1 is straightforward to verify. Notice that  $x_{t+1}=\theta_* z_t +w_t$, and $w_t\mid \F_t=w_t$. By Assumption \ref{ass: w Sigma},  $w_t$ has zero  mean and is bounded by $\|w_t\|_2\leq w_{\max}=\sigma_w \bar w$, thus it is $\sigma_w^2 \bar w^2$-sub-Gaussian, thus satisfying   Condition 1. Condition 2 is verified in Lemma \ref{lem: verify bmsb strictly safe perturbed policy}. 
Condition 3 can be verified below. Notice that
\begin{align*}
	\sigma_{\max}(z_t z_t^\top )\leq \text{trace}(z_t z_t^\top) = \|z_t\|_2^2\leq b_z^2,
\end{align*}
where the last inequality is by Assumption \ref{ass: bounded state action}. Hence, we have 
$\Pb(\sum_{t=1}^T z_t z_t^\top \not \preceq T b_z^2I_{n+m})=0< \delta$ for any $\delta >0$.
Consequently, by applying Theorem \ref{thm: theorem 2.4 general estimation error}, we have
\begin{align*}
	T_0\lesssim&\frac{n+m}{p_z^2}\left(\log(\frac{1}{\delta})+\log(10/p_z)+\log (b_z^2 /s_z^2)\right)\\
	\lesssim &  (m\!+\!n)\max(\bar w^4\!, \bar \eta^4\!)\log\left(\frac{b_z}{\delta}\text{poly}_1(\bar w, \!\bar \eta, \sigma_w^{-1}, \sigma_{\eta}^{-1})) \right),
\end{align*}
where $\text{poly}_1(\bar w, \!\bar \eta, \sigma_w^{-1}, \sigma_{\eta}^{-1}\!)\!=\!\max(\frac{\bar w}{\sigma_w}, \frac{\bar \eta}{\sigma_{\eta}},\frac{\bar w\bar \eta}{\sigma_w\sigma_{\eta}} )\max(\bar w, \bar \eta)$.
The estimation error bound can be organized by:
\begin{align*}
	\|\hat \theta - \theta_*\|_2& \lesssim \frac{\sigma_w \bar w}{ \sqrt T s_z }\sqrt{T_0}\\
	& \lesssim \frac{b_z\sqrt{m+n}}{\sqrt T \sigma_{\eta}} \text{poly}_2(\bar w, \bar \eta, \sigma_w, \sigma_{{\eta}})\\
	& \quad \cdot\log\sqrt{\frac{b_z}{\delta}\text{poly}_1(\bar w, \!\bar \eta, \sigma_w^{-1}, \sigma_{\eta}^{-1})) },
\end{align*}
where $\text{poly}_2(\bar w, \bar \eta, \!\sigma_w,\! \sigma_{{\eta}}\!)\!=\!\bar w \max(\bar w, \!\bar \eta)^2 \max(\bar w\sigma_{\eta}, \bar \eta \sigma_w, \bar w \bar \eta)$. 

\vspace{-3pt}
\subsection{Proof of Lemma \ref{lem: verify bmsb strictly safe perturbed policy}}\label{subsec: verify bmsb}
In this proof, we will first show that the random noises $w_t$ and $\eta_t$ satisfy certain small ball properties, then leverage the properties of $w_t$ and $\eta_t$ to prove the BMSB condition for $\{z_t\}_{t\geq 0}$.

Firstly, we provide the following small-ball properties for $w_t$ and $\eta_t$. 
\begin{lemma}[Supportive lemma]\label{lem: wt etat small ball}
	For any $w_t$ satisfying Assumption \ref{ass: w Sigma}, we have 
	$$ \Pb(\lambda^\top w_t \geq s_w)\geq p_w$$
	for any $\|\lambda\|_2=1$, where $s_w= \frac{\sigma_w}{4\bar w}$, $p_w=\frac{1}{4\bar w^2}$. 
	
	Similarly, for any $\eta_t$ satisfying Assumption \ref{ass: eta}, we have 
	$\Pb(\lambda^\top \eta_t \geq s_\eta)\geq p_\eta$
	for any $\|\lambda\|_2=1$, where $s_\eta= \frac{\sigma_\eta}{4\bar \eta}$, $p_\eta=\frac{1}{4\bar \eta^2}$. 
\end{lemma}
The proof is deferred to Appendix \ref{subsec: wt, etat small ball}.

Secondly, we leverage the  properties for $w_t$ and $\eta_t$ to prove the BMSB condition for $\{z_t\}_{t\geq 0}$. This is achieved by discussing three cases to be specified below.

\nbf{Preparations. }	For notational simplicity, we define filtrations $\F_t^m=\F(w_0, \dots, w_{t-1},\eta_0, \dots, \eta_{t-1})$. Notice that the policy in Theorem \ref{thm: general estimation error bdd} can be written as $u_t=\pi_t(\F_t^m)+\eta_t$. Remember that $\F_t =\{w_0, \dots, w_{t-1}, \eta_0, \dots, \eta_t\}$, so we have
$z_t \in \F_t$ and
\begin{align*}
	z_{t+1}\mid \F_t=	
	\begin{bmatrix}
		x_{t+1}\\
		u_{t+1}
	\end{bmatrix}\mid \F_t= 
	\begin{bmatrix}
		\theta_* z_t+w_t \mid \F_t\\
		\pi_{t+1}(\F_{t+1}^m)+\eta_{t+1}\mid \F_t
	\end{bmatrix},
\end{align*}

When conditioning on $\F_t$, the vector $\theta_* z_t$ is determined, but the vector $	\pi_{t+1}(\F_{t+1}^m)$ is still random due to the randomness of $w_t$.

For the rest of the proof, we will always condition on $\F_t$. Therefore, we will omit the conditioning notation,  i.e., $\cdot\mid \F_t$,  for notational simplicity.

For notational simplicity, we define $k_0=\max(2/\sqrt 3, 4b_z/s_w)$ and split $\lambda$ by $\lambda=(\lambda_1^\top, \lambda_2^\top)^\top\in \R^{m+n}$, where $\lambda_1\in \R^n$, $\lambda_2\in \R^m$, $\|\lambda\|_2^2=\|\lambda_1\|_2^2+\|\lambda_2\|_2^2=1$. 

 We consider three cases:
 \begin{enumerate}
 	\item[(i)] when $\|\lambda_2\|_2\leq 1/k_0$ and $\lambda_1^\top \theta_* z_t \geq 0$,
 	\item[(ii)] when $\|\lambda_2\|_2\leq 1/k_0$ and $\lambda_1^\top \theta_* z_t < 0$,
 	\item[(iii)] when $\|\lambda_2\|_2> 1/k_0$.
 \end{enumerate} 
In the following, we will show 
$\Pb(|\lambda^\top z_{t+1}| \geq s_z)\geq p_z$ in these three cases, which will complete the proof.

The intuition behind the proof is the following. If $\|\lambda_2\|_2$ is small (Cases 1-2),  the impact of $\lambda_2^\top( \pi_t(\F_t^m)+\eta_t)$ will also be small because $u_t=\pi_t(\F_t^m)+\eta_t$ is bounded, so we can leverage the randomness of $w_t$ to take care of the general policy $\pi_t$. If $\|\lambda_2\|_2$ is large (Case 3), $\lambda_2^\top \eta_t$ is also large, so we can leverage the randomness of $\eta_t$ to take care of the general policy $\pi_t$. The proof details are provided below.

\nbf{Case 1: when $\|\lambda_2\|_2\leq 1/k_0$ and $\lambda_1^\top \theta_* z_t \geq 0$}
\begin{align*}
	\lambda_1^\top w_t& \leq \lambda_1^\top(w_t+\theta_*z_t)\leq | \lambda_1^\top(w_t+\theta_*z_t)|\\
	&=| \lambda^\top z_{t+1}-\lambda_2^\top u_{t+1}|\\
&	\leq | \lambda^\top z_{t+1}|+|\lambda_2^\top u_{t+1}|\\
	&\leq | \lambda^\top z_{t+1}|+ \|\lambda_2\|_2 b_z\\
	& \leq  | \lambda^\top z_{t+1}|+b_z/k_0\leq  | \lambda^\top z_{t+1}|+s_w/4
\end{align*}
where the last inequality uses $k_0\geq 4b_z/s_w$.

Further, notice that $k_0\geq 2/\sqrt 3$, so $\|\lambda_2\|_2^2\leq 1/k_0^2\leq 3/4$, thus, $\|\lambda_1\|_2^2 \geq 1/4$, which means $\|\lambda_1\|_2\geq 1/2$. Therefore,
\begin{align*}
	\Pb(\lambda_1^\top w_t\geq s_w/2)&=\Pb(\frac{\lambda_1^\top w_t}{\|\lambda_1\|_2}\geq \frac{s_w}{2\|\lambda_1\|_2})\\
	&\geq \Pb(\frac{\lambda_1^\top w_t}{\|\lambda_1\|_2}\geq s_w)=p_w
\end{align*}
by Lemma \ref{lem: wt etat small ball}. 

By applying the two inequalities above, we obtain the following.
\begin{align*}
	\Pb(|\lambda^\top z_{t+1}| \geq s_z)&\geq \Pb(|\lambda^\top z_{t+1}| \geq s_w/4)\\
	&= \Pb(|\lambda^\top z_{t+1}|+s_w/4 \geq s_w/2)\\
	&\geq \Pb(\lambda_1^\top w_t\geq s_w/2)\geq p_w
\end{align*}
which completes case 1.

\nbf{Case 2: when $\|\lambda_2\|_2\leq 1/k_0$ and $\lambda_1^\top \theta_* z_t < 0$.} This case can be proved similarly to Case 1. 
\begin{align*}
	\lambda_1^\top w_t& \geq \lambda_1^\top(w_t+\theta_*z_t)\geq -| \lambda_1^\top(w_t+\theta_*z_t)|\\
	&=-| \lambda^\top z_{t+1}-\lambda_2^\top u_{t+1}|\\
	&
	\geq -| \lambda^\top z_{t+1}|-|\lambda_2^\top u_{t+1}|\geq -| \lambda^\top z_{t+1}|- \|\lambda_2\|_2 b_z\\
	& \geq - | \lambda^\top z_{t+1}|-b_z/k_0\geq  -| \lambda^\top z_{t+1}|-s_w/4
\end{align*}
where the last inequality uses $k_0\geq 4b_z/s_w$.

Further, notice that $k_0\geq 2/\sqrt 3$, so $\|\lambda_2\|_2^2\leq 1/k_0^2\leq 3/4$, thus, $\|\lambda_1\|_2^2 \geq 1/4$, which means $\|\lambda_1\|_2\geq 1/2$. Therefore,
\begin{align*}
&\ 	\Pb(\lambda_1^\top w_t\leq -s_w/2)=\Pb(\frac{\lambda_1^\top w_t}{\|\lambda_1\|_2}\leq -\frac{s_w}{2\|\lambda_1\|_2}) \\
	\geq &\  \Pb(\frac{\lambda_1^\top w_t}{\|\lambda_1\|_2}\leq -s_w)=\Pb(\frac{-\lambda_1^\top w_t}{\|\lambda_1\|_2}\geq s_w)=p_w
\end{align*}
by $s_w/(2\|\lambda_1\|_2)\leq s_w$, and thus $-s_w/(2\|\lambda_1\|_2)\geq -s_w$, and Assumption \ref{ass: w Sigma}.

Consequently,
\begin{align*}
	\Pb(|\lambda^\top z_{t+1}| \geq s_z)&\geq \Pb(|\lambda^\top z_{t+1}| \geq s_w/4)\\
	&= \Pb(-|\lambda^\top z_{t+1}|-s_w/4 \leq -s_w/2)\\
	&\geq \Pb(\lambda_1^\top w_t\leq -s_w/2)\geq p_w
\end{align*}
which completes the proof of Case 2.

\nbf{Case 3: when $\|\lambda_2\|_2> 1/k_0$.} Define
\begin{align*}
	\Omega_1^\lambda&=\{w_t\in \R^n\mid \lambda_1^\top(w_t+\theta_*z_t)+\lambda_2^\top(\pi_{t+1}(\F_{t+1}^m))\geq 0\}\\
	\Omega_2^\lambda&=\{w_t\in \R^n\mid \lambda_1^\top(w_t+\theta_*z_t)+\lambda_2^\top(\pi_{t+1}(\F_{t+1}^m))< 0\}
\end{align*}
Notice that $\Pb( w_t \in \Omega_1^\lambda)\!+\!\Pb( w_t \in \Omega_2^\lambda)=1$. Further, we have
\begin{align*}
	\Pb(&|\lambda^\top z_{t+1}| \geq s_z)\geq \Pb(|\lambda^\top z_{t+1}| \geq v)\\
	&= \Pb(\lambda^\top z_{t+1} \geq v)+ \Pb(\lambda^\top z_{t+1} \leq -v)\\
	& \geq \Pb(\lambda^\top z_{t+1} \geq v, w_t \in \Omega_1^\lambda)+ \Pb(\lambda^\top z_{t+1} \leq -v, w_t \in \Omega_2^\lambda)\\
	& \geq \Pb(\lambda_2^\top \eta_{t+1} \geq v, w_t \in \Omega_1^\lambda)+\Pb(\lambda_2^\top \eta_{t+1} \leq -v, w_t \in \Omega_2^\lambda)\\
	&=\Pb(\lambda_2^\top \eta_{t+1} \geq v)\Pb( w_t \in \Omega_1^\lambda)\\
	& \quad +\Pb(\lambda_2^\top \eta_{t+1} \leq -v)\Pb( w_t \in \Omega_2^\lambda)\geq p_\eta,
\end{align*}
where $v= s_\eta/k_0=\min(\sqrt 3  s_\eta/2, s_w  s_\eta/(4 b_z))$ and the last inequality is because of the following arguments.
Notice that, by Lemma \ref{lem: wt etat small ball}, we have
\begin{align*}
	\Pb(\lambda_2^\top \eta_{t+1} \geq v)&= \Pb(\lambda_2^\top \eta_{t+1}/\|\lambda_{2}\|_2 \geq v/\|\lambda_{2}\|_2)\\
	& \geq \Pb(\lambda_2^\top  \eta_{t+1}/\|\lambda_{2}\|_2 \geq k_0v)\\
	& = \Pb(\lambda_2^\top  \eta_{t+1}/\|\lambda_{2}\|_2 \geq  s_\eta)\geq p_\eta.
\end{align*}
Similarly, we can obtain
$\Pb(\lambda_2^\top \eta_{t+1} \leq -v)=	\Pb(-\lambda_2^\top \eta_{t+1} \geq v)\geq p_\eta.$
This completes the proof of Case 3.

By combining Cases 1-3, we completed the proof. \qed


%% file: sec4_theory_proof.tex
\section{Applications to safe learning for  constrained LQR}\label{sec: rmpc}
This section will introduce the applications of our system identification error bound to safe learning of constrained LQR. In particular, we will use RMPC as an illustrative example. Other safe control policies reviewed in Section \ref{sec: problem} can be applied similarly. 

Firstly, we introduce a constrained LQR problem with model uncertainties below. 
\begin{equation}\label{equ: J(pi)}
	\begin{aligned}
		&	\min_{u_0, u_1, \dots}\lim_{T\to +\infty}\frac{1}{T}\sum_{t=0}^{T}\E [x_t^\top Q x_t+ u_t^\top R u_t]\\
		\text{s.t.} &\ \ x_{t+1}=A_*x_t+B_* u_t +w_t, \ \forall \, t\geq 0,\\
	&\ 	\  x_t \in \X, \ u_t \in \U, \    \forall  \, t\geq 0,\   \forall\, \{w_k \in \W\}_{k\geq 0}.
	\end{aligned}
\end{equation}
where  the system parameters $\theta_*= (A_*, B_*)$ are not accurately known. However, in the robust control framework, some domain knowledge of $\theta_*$ is usually assumed to be known. In particular, a bounded uncertainty set $\Theta_0$ is usually assumed to be known and to contain the true parameters, i.e., $\theta_*\in \Theta_0$ \citep{lorenzen2019robust,kohler2019linear,rawlings2009model}.

Next, we briefly review RMPC below. For simplicity, we will only introduce the basic form of tube-based RMPC in \cite{rawlings2009model} below and note that there have been significant efforts on improving the basic form, e.g., \citep{lorenzen2019robust,kohler2019linear}. The high-level intuition behind tube-based RMPC is to plan a nominal trajectory, denoted by $x_{t+k|t}$, and construct a tube $\Sb_K$ such that the true trajectory  $x_{t+k}$ always lies within the tube around the nominal trajectory, i.e., $x_{t+k}\in x_{t+k|k}\oplus \Sb_K$ (constraints on $u_t$ are handled similarly). Then, by requiring the tube around the nominal trajectory to satisfy the constraints, tube-based RMPC achieves robust constraint satisfaction despite uncertainties in the system. 
In particular, RMPC solves the following finite-horizon optimal control problem to obtain $\{v_{t|t}^*, \dots, v_{t+W-1|t}^*\}$ and implements the control input $u_t=\pi_{\text{RMPC}}(x_t)=Kx_t + v_{t|t}^*$ at each time $t$, where $K$ is introduced below.
\begin{align}
		&	\min_{\!\{\!v_{t\!+\!k|t}\!\}_{k=0}^{W\!-\!1}}\!\sum_{k=0}^{W\!-\!1}\!\E [x_{t+k|t}^\top Q x_{t+\!k|t}\!+ \! u_{t+\!k|t}^\top R u_{t+\!k|t}]\!\!+\!\! V_f(x_{t+W|t}\!)\notag\\
		&\begin{aligned}\label{equ: rmpc}
		\text{s.t.  }  &x_{t+k+1|t}=A_0 x_{t+k|t}+B_0 u_{t+k|t},  \forall \, 0\leq k\leq W\!-\!1\\
		&   u_{t+k|t} = K x_{t+k|t}+ v_{t+k|t}, \ \forall \, 0\leq k \leq W\!-\!1\\
		& x_{t+k|t}\! \in \!\X\!\ominus \!\Sb_K, \, u_{t+k|t} \!\in \!\U\!\ominus\! K\Sb_K,  \forall \, 0\!\leq\! k\leq \!W\!-\!1\\
		&	 x_{t|t}=x_t, \ x_{t+W|t}\in \X_f\subseteq \in \X\ominus \Sb_K
		\end{aligned}
	\end{align}
where  the feedback gain $K$ is assumed to stabilize all the systems in $\Theta_0$, the initial system estimation is $(A_0, B_0)\in \Theta_0$, the terminal cost  $V_f(\cdot)$ and terminal constraint $\X_f$ needs to satisfy the assumptions in Section 3.5 of \cite{rawlings2009model}, and the tube $\Sb_K$ is defined as
\begin{equation}\label{equ: tube}
\begin{aligned}
	\Sb_K & = \sum_{i=0}^{+\infty} (A_0+B_0K)^i \Sb\\
		\Sb&=\{w+(\theta-\theta_0)z: w\in \W, x\in \X, u\in \U, \theta\in \Theta\}
\end{aligned}
\end{equation}
It is worth mentioning that the tube design in \eqref{equ: tube} is  conservative and can be  improved in more advanced RMPC methods, e.g., \cite{lorenzen2019robust,kohler2019linear}.


To learn the true parameters $\theta_*$, we introduce random noises $\eta_t \in\mathbb  H=\{\eta: \|\eta\|_2\leq \eta_{\max}\}$ to provide enough excitation. In particular, we consider policy $u_t=\pi_{\text{RMPC}}(x_t)+\eta_t$. Due to the additional noises, we need to adjust the tube to account for the additional uncertainty by the following.
\begin{align}\label{equ: eta tube}
	\Sb_{\eta}&=\Sb \oplus \mathbb H,\ \ 	\Sb_{K, \eta} = \sum_{i=0}^{+\infty} (A_0+B_0K)^i \Sb_{\eta}.
\end{align}
Then, we can retain the robust constraint satisfaction of policy $u_t=\pi_{\text{RMPC}}(x_t)+\eta_t$ and apply our Theorem \ref{thm: theorem 2.4 general estimation error}.

Further,  we can  adjust our LSE estimation  to be consistent with the prior knowledge that $\theta_*\in \Theta_0$. In particular, we can obtain a point estimator $\tilde \theta =\argmin_{\theta \in \Theta_0}\|\hat \theta-\theta\|_2^2$ by projection with the same  estimation error bound. The details are provided in the corollary below.
\begin{corollary}
	Consider a single trajectory $\{x_t, u_t\}_{t=0}^T$ generated by RMPC with excitation $u_t=\pi_{\text{RMPC}}(x_t)+\eta_t$, where the tube $\Sb_K$ in \eqref{equ: rmpc} are adjusted to $\Sb_{K,\eta}$ in \eqref{equ: eta tube}. Suppose the constraints $\X, \U$ are bounded, i.e., there exists $b_z=\max_{x\in \X, u\in \U}\sqrt{\|x\|_2^2+ \|u\|^2} $. If Assumptions \ref{ass: w Sigma} and \ref{ass: eta}  are true and  the condition on $T$ in Theorem \ref{thm: theorem 2.4 general estimation error} is satisfied, then 
$
		\|\tilde \theta - \theta_*\|_2 \lesssim \frac{b_z\sqrt{m+n}}{\sqrt T \sigma_{\eta}} \text{poly}_2(\bar w, \bar \eta, \sigma_w, \sigma_{{\eta}})\\
	\times	\sqrt{\log(\frac{b_z}{\delta})\!+\!\log(\text{poly}_1(\bar w, \!\bar \eta, \sigma_w^{-1}, \sigma_{\eta}^{-1})) )}.$
%

\end{corollary}
\textit{Proof: } The boundedness of states and actions follows directly from the constraint satisfaction of RMPC. Further, since $\theta_*\in \Theta_0$, by the non-expansiveness of projection, we have $\|\tilde \theta - \theta_*\|_2\leq \|\hat \theta - \theta_*\|_2$. Then, the proof is completed by applying Theorem \ref{thm: general estimation error bdd}. \qed

It is worth mentioning that for computational purposes, the projection with respect to the Frobenius norm can also be adopted here. In this case, the estimation error bound will increase by a factor $\sqrt n$ due to the change of norms.

%% file: sec5_numerical.tex
\section{Numerical experiments}

%

\begin{figure}
	\centering
	\begin{subfigure}{.25\textwidth}
		\centering
			\includegraphics[width=.98\linewidth]{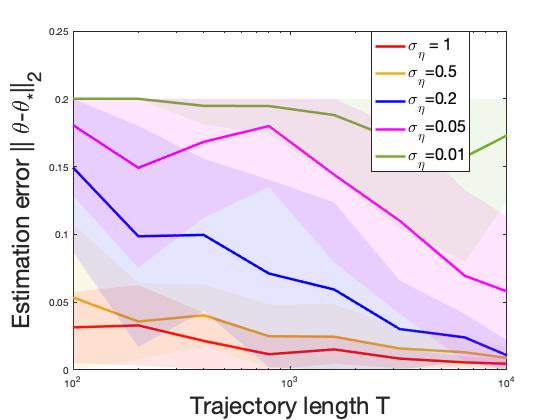}
		\caption{Constrained LQR}
	\end{subfigure}%
	\begin{subfigure}{.25\textwidth}
		\centering

			\includegraphics[width=.98\linewidth]{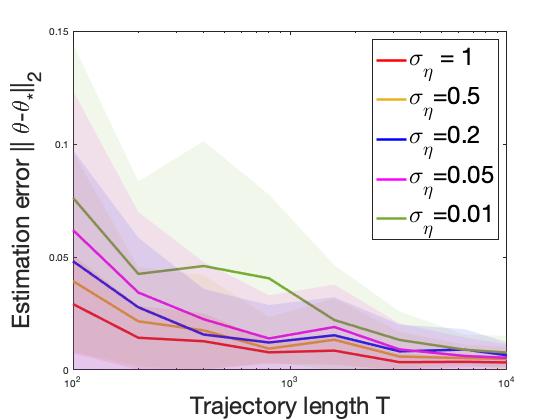}
		\caption{Constrained LQ Tracking}
	\end{subfigure}
	\caption{The figures  plot the estimation error of $\tilde \theta$ when applying  $u_t = \pi_{\text{RMPC}}(x_t)+\eta_t$ under different excitation levels   $\sigma_{\eta}$. Figure (a) considers the constrained LQR  in  \eqref{equ: J(pi)}, so $\pi_{\text{RMPC}}(x_t)$ is time-invariant. Figure (b) considers a time-varying tracking problem, so $\pi_{\text{RMPC}}(x_t)$ is time-varying. The solid lines represent the sample mean. The shades represent one standard deviation.}
	\label{fig:estimation error}
\end{figure}

\begin{figure}
	\centering
	\begin{subfigure}{.25\textwidth}
		\centering
		\includegraphics[width=.98\linewidth]{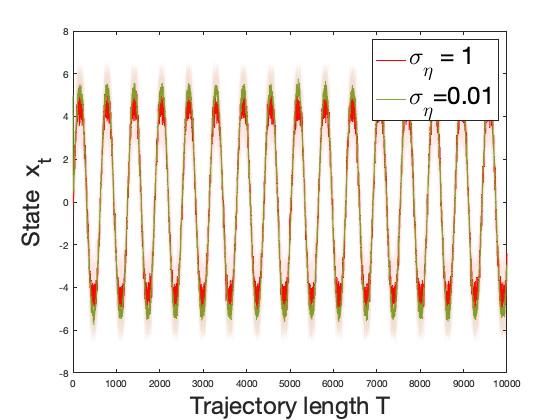}
		\caption{$x_t$ trajectories}
	\end{subfigure}%
	\begin{subfigure}{.25\textwidth}
		\centering
		\includegraphics[width=.98\linewidth]{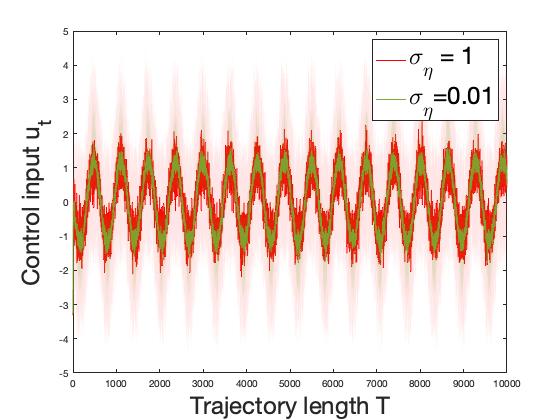}
		\caption{$u_t$ trajectories}
	\end{subfigure}
	\caption{The solid lines represent the sample means and the shades contain all possible trajectories of states and actions generated by RMPC in the tracking problem.}
	\label{fig:trajectory}
\end{figure}

In this section, we provide numerical experiments to supplement our theoretical analysis by learning with RMPC policies as reviewed in Section \ref{sec: rmpc}.

In our experiment, we consider a linear system \eqref{equ: LTI system} with $A_*=1.2$, $B_*=0.9$, and a model uncertainty set $\Theta_0=[1,1.2]\times [0.9,1.1]$. The system disturbances $w_t$ are i.i.d. following a Uniform distribution on $[-1,1]$. We apply the  basic tube-based RMPC policy reviewed in Section \ref{sec: rmpc} and  refer the reader to \cite{rawlings2009model} for more details. We let the initial estimator in RMPC \eqref{equ: rmpc} be $A_0=1.1, B_0=1$. We generate the excitation noises by $\eta_t = \sigma_{\eta} \tilde \eta_t$, where $\sigma_{\eta}$ represents the excitation level and the variance of $\eta_t$, and $\tilde \eta_t$ are i.i.d. generated from a Uniform distribution on $\{-1,1\}$. We consider the state constraint and the control input constraints as $[-10, 10]$. We let the RMPC lookahead window be $W=5$. We repeat each experiment for 15 times.

We consider two types of problems: (a) the constrained LQR problem as reviewed in Section \ref{sec: rmpc} and (b) the constrained LQ tracking problem with a time-varying cost function  $(x_t-g_t)^\top Q (x_t-g_t)+u_t^\top R u_t$, where we generate the target trajectory by $g_t=8 \sin(t/100)$.  In the constrained LQR problem, the RMPC policy is time-invariant and nonlinear. In the constrained LQ tracking problem, the RMPC policy is both time-varying and nonlinear.

In Figure \ref{fig:estimation error}, we plot the estimation errors of both constrained LQR and constrained LQ tracking under different excitation levels.  We  observe that,  in both cases, the estimation errors decrease with the trajectory lengths $T$. Besides, the estimation errors tend to be smaller if the excitation level $\sigma_{\eta}$ is larger. Both observations above are consistent with Theorem \ref{thm: general estimation error bdd}. It is worth noting that the excitation level $\sigma_{\eta}$ cannot be too large otherwise,   the RMPC problem \eqref{equ: rmpc}  becomes infeasible. Interestingly, in this case, the LQ tracking problem yields a smaller estimation error. This is because the tracking of the moving target helps with the system exploration  in this setting.

In Figure \ref{fig:trajectory}, we plot the state and control input trajectories in the constrained LQ tracking problem under  different levels of excitation noises.  Figure \ref{fig:trajectory} demonstrates that RMPC guarantees constraint satisfaction under the model uncertainties and excitation noises even when the target   drives the state  towards the boundaries of the constraints, thus validating our Assumption \ref{ass: bounded state action}. Further, with a larger excitation level $\sigma_{\eta}$, the possible region of the trajectories is larger due to more uncertainties in the system.

%% file: sec6_conclusion.tex
\section{Conclusion}
This paper studies the linear system identification  by a single-trajectory of data generated by general nonlinear and/or time-varying policies with i.i.d. random excitations. We provide a general estimation error bound for any policies with bounded  states and actions. Our bound for general policies is consistent with that for linear policies with respect to the trajectory length, system dimensions, and excitation levels. We apply our results to safe learning with robust model predictive control and conduct numerical experiments. There are many future directions to explore, e.g., (i) applying our results to the adaptive learning of robust MPC to determine the tube sizes and conduct regret analysis, (ii) relaxing the bounded disturbances and bounded trajectories assumptions to (sub)Gaussian disturbances and chance constraints of trajectories, (iii) understanding the fundamental lower bound of this problem, (iv) exploring what structures of nonlinear systems can provide similar identification guarantees, (v) designing active exploration with better estimation performance, and (iv) studying the estimation guarantees of other methods, e.g., set membership.

%
%
%
%
%
%
%
%

\appendix
\section{Proof of Lemma \ref{lem: wt etat small ball}}\label{subsec: wt, etat small ball}
We first consider $w_t$. For any fixed $\lambda$ such that $\|\lambda\|_2=1$, we define $y=\lambda^\top w_t$. By Assumption \ref{ass: w Sigma}, we have
$
	\E y^2=E \lambda^\top w_t w_t^\top \lambda = \lambda^\top \text{cov}(w_t) \lambda \geq \sigma_w^2
$
and $|y|\leq \|\lambda\|_2\|w_t\|_2 \leq w_{\max}$. Therefore, 
$\E |y| \geq \E y^2/w_{\max} \geq \sigma_w^2/w_{\max}.$
By leveraging the inequality above and $\E y=0$, we obtain
$
	\E y \one_{(y\geq 0)} = \frac{1}{2}(\E|y| +\E y) \geq \frac{\sigma_w^2}{2w_{\max}}
$.
Further, we have
\begin{align*}
\frac{\sigma_w^2}{2w_{\max}} &\leq 	\E y \one_{(y\geq 0)}\\
&\leq w_{\max} \Pb(y\geq \frac{\sigma_w^2}{4w_{\max}})+ \frac{\sigma_w^2}{4w_{\max}} \Pb(0\leq y< \frac{\sigma_w^2}{4w_{\max}})\\
	& \leq w_{\max} \Pb(y\geq \frac{\sigma_w^2}{4w_{\max}})+ \frac{\sigma_w^2}{4w_{\max}} 
\end{align*}
By rearranging the terms, we obtain $\Pb(y \geq s_w)\geq p_w$, where $s_w= \frac{\sigma_w^2}{4w_{\max}}=\frac{\sigma_w}{4\bar w}$ and $p_w=\frac{\sigma_w^2}{4w_{\max}^2}=\frac{1}{4\bar w^2}$. The proof for $\eta_t$ is the same.\qed

\vspace{-3pt}

%% file: main.bbl
\begin{thebibliography}{32}
\providecommand{\natexlab}[1]{#1}
\providecommand{\url}[1]{\texttt{#1}}
\providecommand{\urlprefix}{URL }
\expandafter\ifx\csname urlstyle\endcsname\relax
  \providecommand{\doi}[1]{doi:\discretionary{}{}{}#1}\else
  \providecommand{\doi}{doi:\discretionary{}{}{}\begingroup
  \urlstyle{rm}\Url}\fi

\bibitem[{Bai et~al.(1998)Bai, Cho, and Tempo}]{bai1998convergence}
Bai, E.W., Cho, H., and Tempo, R. (1998).
\newblock Convergence properties of the membership set.
\newblock \emph{Automatica}, 34(10), 1245--1249.

\bibitem[{Boyd and Sastry(1986)}]{boyd1986necessary}
Boyd, S. and Sastry, S.S. (1986).
\newblock Necessary and sufficient conditions for parameter convergence in
  adaptive control.
\newblock \emph{Automatica}, 22(6), 629--639.

\bibitem[{Dean et~al.(2018)Dean, Mania, Matni, Recht, and Tu}]{dean2018regret}
Dean, S., Mania, H., Matni, N., Recht, B., and Tu, S. (2018).
\newblock Regret bounds for robust adaptive control of the linear quadratic
  regulator.
\newblock In \emph{Advances in Neural Information Processing Systems},
  4188--4197.

\bibitem[{Dean et~al.(2019{\natexlab{a}})Dean, Mania, Matni, Recht, and
  Tu}]{dean2019sample}
Dean, S., Mania, H., Matni, N., Recht, B., and Tu, S. (2019{\natexlab{a}}).
\newblock On the sample complexity of the linear quadratic regulator.
\newblock \emph{Foundations of Computational Mathematics}, 1--47.

\bibitem[{Dean et~al.(2019{\natexlab{b}})Dean, Tu, Matni, and
  Recht}]{dean2019safely}
Dean, S., Tu, S., Matni, N., and Recht, B. (2019{\natexlab{b}}).
\newblock Safely learning to control the constrained linear quadratic
  regulator.
\newblock In \emph{2019 American Control Conference (ACC)}, 5582--5588. IEEE.

\bibitem[{Dogan et~al.(2021)Dogan, Shen, and Aswani}]{dogan2021regret}
Dogan, I., Shen, Z.J.M., and Aswani, A. (2021).
\newblock Regret analysis of learning-based mpc with partially-unknown cost
  function.
\newblock \emph{arXiv preprint arXiv:2108.02307}.

\bibitem[{Fisac et~al.(2018)Fisac, Akametalu, Zeilinger, Kaynama, Gillula, and
  Tomlin}]{fisac2018general}
Fisac, J.F., Akametalu, A.K., Zeilinger, M.N., Kaynama, S., Gillula, J., and
  Tomlin, C.J. (2018).
\newblock A general safety framework for learning-based control in uncertain
  robotic systems.
\newblock \emph{IEEE Transactions on Automatic Control}, 64(7), 2737--2752.

\bibitem[{Fogel and Huang(1982)}]{fogel1982value}
Fogel, E. and Huang, Y.F. (1982).
\newblock On the value of information in system identification—bounded noise
  case.
\newblock \emph{Automatica}, 18(2), 229--238.

\bibitem[{Foster et~al.(2020)Foster, Sarkar, and Rakhlin}]{foster2020learning}
Foster, D., Sarkar, T., and Rakhlin, A. (2020).
\newblock Learning nonlinear dynamical systems from a single trajectory.
\newblock In \emph{Learning for Dynamics and Control}, 851--861. PMLR.

\bibitem[{K{\"o}hler et~al.(2019)K{\"o}hler, Andina, Soloperto, M{\"u}ller, and
  Allg{\"o}wer}]{kohler2019linear}
K{\"o}hler, J., Andina, E., Soloperto, R., M{\"u}ller, M.A., and Allg{\"o}wer,
  F. (2019).
\newblock Linear robust adaptive model predictive control: Computational
  complexity and conservatism.
\newblock In \emph{2019 IEEE 58th Conference on Decision and Control (CDC)},
  1383--1388. IEEE.

\bibitem[{Li et~al.(2021{\natexlab{a}})Li, Das, and Li}]{li2021online}
Li, Y., Das, S., and Li, N. (2021{\natexlab{a}}).
\newblock Online optimal control with affine constraints.
\newblock In \emph{Proceedings of the AAAI Conference on Artificial
  Intelligence}, volume~35, 8527--8537.

\bibitem[{Li et~al.(2021{\natexlab{b}})Li, Das, Shamma, and Li}]{li2021safe}
Li, Y., Das, S., Shamma, J., and Li, N. (2021{\natexlab{b}}).
\newblock Safe adaptive learning-based control for constrained linear quadratic
  regulators with regret guarantees.
\newblock \emph{arXiv preprint arXiv:2111.00411}.

\bibitem[{Lim{\'o}n et~al.(2010)Lim{\'o}n, Alvarado, Alamo, and
  Camacho}]{limon2010robust}
Lim{\'o}n, D., Alvarado, I., Alamo, T., and Camacho, E.F. (2010).
\newblock Robust tube-based mpc for tracking of constrained linear systems with
  additive disturbances.
\newblock \emph{Journal of Process Control}, 20(3), 248--260.

\bibitem[{Lopez et~al.(2020)Lopez, Slotine, and How}]{lopez2020robust}
Lopez, B.T., Slotine, J.J.E., and How, J.P. (2020).
\newblock Robust adaptive control barrier functions: An adaptive and
  data-driven approach to safety.
\newblock \emph{IEEE Control Systems Letters}, 5(3), 1031--1036.

\bibitem[{Lorenzen et~al.(2019)Lorenzen, Cannon, and
  Allg{\"o}wer}]{lorenzen2019robust}
Lorenzen, M., Cannon, M., and Allg{\"o}wer, F. (2019).
\newblock Robust mpc with recursive model update.
\newblock \emph{Automatica}, 103, 461--471.

\bibitem[{Mania et~al.(2022)Mania, Jordan, and Recht}]{mania2022active}
Mania, H., Jordan, M.I., and Recht, B. (2022).
\newblock Active learning for nonlinear system identification with guarantees.
\newblock \emph{Journal of Machine Learning Research}, 23(32), 1--30.

\bibitem[{Mhammedi et~al.(2020)Mhammedi, Foster, Simchowitz, Misra, Sun,
  Krishnamurthy, Rakhlin, and Langford}]{mhammedi2020learning}
Mhammedi, Z., Foster, D.J., Simchowitz, M., Misra, D., Sun, W., Krishnamurthy,
  A., Rakhlin, A., and Langford, J. (2020).
\newblock Learning the linear quadratic regulator from nonlinear observations.
\newblock \emph{Advances in Neural Information Processing Systems}, 33,
  14532--14543.

\bibitem[{Oldewurtel et~al.(2008)Oldewurtel, Jones, and
  Morari}]{oldewurtel2008tractable}
Oldewurtel, F., Jones, C.N., and Morari, M. (2008).
\newblock A tractable approximation of chance constrained stochastic mpc based
  on affine disturbance feedback.
\newblock In \emph{2008 47th IEEE conference on decision and control},
  4731--4736. IEEE.

\bibitem[{Oymak(2019)}]{oymak19stochastic}
Oymak, S. (2019).
\newblock Stochastic gradient descent learns state equations with nonlinear
  activations.
\newblock In A.~Beygelzimer and D.~Hsu (eds.), \emph{Proceedings of the
  Thirty-Second Conference on Learning Theory}, volume~99 of \emph{Proceedings
  of Machine Learning Research}, 2551--2579. PMLR.

\bibitem[{Oymak and Ozay(2019)}]{oymak2019non}
Oymak, S. and Ozay, N. (2019).
\newblock Non-asymptotic identification of lti systems from a single
  trajectory.
\newblock In \emph{2019 American control conference (ACC)}, 5655--5661. IEEE.

\bibitem[{Rawlings and Mayne(2009)}]{rawlings2009model}
Rawlings, J.B. and Mayne, D.Q. (2009).
\newblock \emph{Model predictive control: Theory and design}.
\newblock Nob Hill Pub.

\bibitem[{Salehi et~al.(2022)Salehi, Taplin, and Dani}]{salehi2022learning}
Salehi, I., Taplin, T., and Dani, A. (2022).
\newblock Learning discrete-time uncertain nonlinear systems with probabilistic
  safety and stability constraints.
\newblock \emph{IEEE Open Journal of Control Systems}.

\bibitem[{Sarkar and Rakhlin(2019)}]{sarkar2019near}
Sarkar, T. and Rakhlin, A. (2019).
\newblock Near optimal finite time identification of arbitrary linear dynamical
  systems.
\newblock In \emph{International Conference on Machine Learning}, 5610--5618.
  PMLR.

\bibitem[{Sattar and Oymak(2022)}]{sattar2022non}
Sattar, Y. and Oymak, S. (2022).
\newblock Non-asymptotic and accurate learning of nonlinear dynamical systems.
\newblock \emph{Journal of Machine Learning Research}, 23(140), 1--49.

\bibitem[{Sattar et~al.(2022)Sattar, Oymak, and Ozay}]{sattar2022finite}
Sattar, Y., Oymak, S., and Ozay, N. (2022).
\newblock Finite sample identification of bilinear dynamical systems.
\newblock \emph{arXiv preprint arXiv:2208.13915}.

\bibitem[{Simchowitz et~al.(2018)Simchowitz, Mania, Tu, Jordan, and
  Recht}]{simchowitz2018learning}
Simchowitz, M., Mania, H., Tu, S., Jordan, M.I., and Recht, B. (2018).
\newblock Learning without mixing: Towards a sharp analysis of linear system
  identification.
\newblock In \emph{Conference On Learning Theory}, 439--473. PMLR.

\bibitem[{Taylor et~al.(2020)Taylor, Singletary, Yue, and
  Ames}]{taylor2020control}
Taylor, A.J., Singletary, A., Yue, Y., and Ames, A.D. (2020).
\newblock A control barrier perspective on episodic learning via
  projection-to-state safety.
\newblock \emph{IEEE Control Systems Letters}, 5(3), 1019--1024.

\bibitem[{Wabersich and Zeilinger(2018)}]{wabersich2018linear}
Wabersich, K.P. and Zeilinger, M.N. (2018).
\newblock Linear model predictive safety certification for learning-based
  control.
\newblock In \emph{2018 IEEE Conference on Decision and Control (CDC)},
  7130--7135. IEEE.

\bibitem[{Wabersich and Zeilinger(2020)}]{wabersich2020performance}
Wabersich, K.P. and Zeilinger, M.N. (2020).
\newblock Performance and safety of bayesian model predictive control: Scalable
  model-based rl with guarantees.
\newblock \emph{arXiv preprint arXiv:2006.03483}.

\bibitem[{Xu(2018)}]{xu2018constrained}
Xu, X. (2018).
\newblock Constrained control of input--output linearizable systems using
  control sharing barrier functions.
\newblock \emph{Automatica}, 87, 195--201.

\bibitem[{Zhao and Li(2022)}]{zhao2022adaptive}
Zhao, Z. and Li, Q. (2022).
\newblock Adaptive sampling methods for learning dynamical systems.
\newblock In \emph{Mathematical and Scientific Machine Learning}, 335--350.
  PMLR.

\bibitem[{Ziemann and Tu(2022)}]{ziemann2022learning}
Ziemann, I. and Tu, S. (2022).
\newblock Learning with little mixing.
\newblock \emph{arXiv preprint arXiv:2206.08269}.

\end{thebibliography}
